\documentclass[reqno,12pt]{amsart}
\usepackage{amsmath, amsthm, amssymb, booktabs}
\usepackage{mathrsfs}
\usepackage{color}

\advance \topmargin by -\headheight
\advance \topmargin by -\headsep
     
\evensidemargin \oddsidemargin
\newtheorem{theorem}{Theorem}[section]
\newtheorem{lemma}[theorem]{Lemma}

\theoremstyle{definition}

\theoremstyle{remark}

\numberwithin{equation}{section}

\newcommand{\mmod}[1]{\,\,(\text{mod}\,\,#1)}

\def\calC{{\mathscr C}}

\def\dbC{{\mathbb C}}\def\dbN{{\mathbb N}}
\def\dbR{{\mathbb R}}
\def\dbZ{{\mathbb Z}}

\def\grB{{\mathfrak B}}

\def\grJ{{\mathfrak J}}
\def\grk{{\mathfrak k}}

\def\grm{{\mathfrak m}}\def\grM{{\mathfrak M}}\def\grN{{\mathfrak N}}
\def\grS{{\mathfrak S}}
\def\grB{{\mathfrak B}}
\def\grK{{\mathfrak K}}

\def\alp{{\alpha}} 
\def\bet{{\beta}}  
\def\gam{{\gamma}} \def\Gam{{\Gamma}}
\def\del{{\delta}} \def\Del{{\Delta}}

\def\tet{{\theta}}  \def\Tet{{\Theta}}
\def\kap{{\kappa}}
\def\lam{{\lambda}}

\def\Ups{{\Upsilon}} 

\def\ome{{\omega}}  
\def\d{{\partial}}
\def\eps{\varepsilon}

\def\le{\leqslant} \def\ge{\geqslant}

\def\d{{\,{\rm d}}}

\makeatletter
\@namedef{subjclassname@2020}{\textup{2020} Mathematics Subject Classification}
\makeatother

\begin{document}
\title[Waring's problem]{On Waring's problem for larger powers}
\author[J\"org Br\"udern]{J\"org Br\"udern}
\address{Mathematisches Institut, Bunsenstrasse 3--5, D-37073 G\"ottingen, Germany}
\email{jbruede@gwdg.de}
\author[Trevor D. Wooley]{Trevor D. Wooley}
\address{Department of Mathematics, Purdue University, 150 N. University Street, West 
Lafayette, IN 47907-2067, USA}
\email{twooley@purdue.edu}
\thanks{First author supported by Deutsche Forschungsgemeinschaft Project Number 
255083470. Second author supported by NSF grants DMS-1854398 and DMS-2001549.}
\date{}

\begin{abstract} Let $G(k)$ denote the least number $s$ having the property that every 
sufficiently large natural number is the sum of at most $s$ positive integral $k$-th powers. 
Then for all $k\in \mathbb N$, one has
\[
G(k)\le \lceil k(\log k+4.20032)\rceil .
\]
Our new methods improve on all bounds available hitherto when $k\ge 14$.
\end{abstract}
\maketitle

\section{Introduction} Since the introduction by Hardy and Littlewood of their circle method 
a century ago (see \cite{HL1920}), it has been possible to surmise progress associated with 
this technology from corresponding advances in the theory of Waring's problem. As is usual, 
we denote by $G(k)$ the least number $s$ having the property that every sufficiently large 
natural number is the sum of at most $s$ positive integral $k$-th powers. The initial bound 
$G(k)\le (k-2)2^{k-1}+5$ of Hardy and Littlewood \cite{HL1922} was improved rapidly over 
the next four decades, culminating in 1959 with Vinogradov's bound
\[
G(k)\le k(2\log k+4\log \log k+2\log \log \log k +13)\quad (k\ge 170,000)
\]
(see \cite{Vin1959}). The latter bound was subsequently improved by Karatsuba 
\cite{Kar1985}, and shortly thereafter by Vaughan \cite{Vau1989}, showing that
\[
G(k)\le 2k(\log k+\log \log k+1+\log 2+O(\log \log k/\log k)).
\]
A little over three decades after the work of Vinogradov, the second author obtained a 
bound roughly half that of this earlier work, establishing the bound
\[
G(k)\le k(\log k+\log \log k+2+O(\log \log k/\log k))
\]
(see \cite{Woo1990, Woo1992} and \cite[Theorem 1.4]{Woo1995}). Our primary goal in 
this memoir is the removal of the secondary term of size $k\log \log k$.

\begin{theorem}\label{theorem1.1} For all $k\in \mathbb N$, one has 
$G(k)\le \lceil k(\log k+4.20032)\rceil $.
\end{theorem}

The conclusion of this theorem constitutes the largest improvement in available bounds for 
$G(k)$, when $k$ is large, since the progress achieved thirty years ago by the second 
author \cite{Woo1990, Woo1992}. The upper bound presented in Theorem 
\ref{theorem1.1} is in fact an approximation to one asymptotically very slightly stronger. In 
order to describe this result, we introduce some auxiliary constants. Let $\ome$ be the 
unique real solution, with $\ome \ge 1$, of the transcendental equation
\begin{equation}\label{1.1}
\ome-2-1/\ome =\log \ome .
\end{equation}
We then put
\begin{equation}\label{1.2}
C_1=2+\log (\ome^2-3-2/\ome)\quad\text{and}\quad 
C_2=\frac{\ome^2+3\ome-2}{\ome^2-\ome-2}.
\end{equation}
A modest computation reveals that
\[
\ome =3.548292\ldots ,\quad C_1=4.200189\ldots \quad \text{and}\quad 
C_2=3.015478\ldots .
\]

\begin{theorem}\label{theorem1.2}
For all $k\in \mathbb N$, one has $G(k)<k(\log k+C_1)+C_2$.
\end{theorem}  

It transpires that the new ideas underlying the progress exhibited in Theorems 
\ref{theorem1.1} and \ref{theorem1.2} apply not only for very large values of $k$, but also 
for exponents of moderate size.

\begin{theorem}\label{theorem1.3}
When $14\le k\le 20$, one has $G(k)\le H(k)$, where $H(k)$ is defined by means of Table 
\ref{tab1}.
\end{theorem}

\begin{table}[h]
\begin{center}
\begin{tabular}{ccccccccccccccccc}
\toprule
$k$ & $14$ & $15$ & $16$ & $17$ & $18$ & $19$ & $20$\\
$H(k)$ & $89$ & $97$ & $105$ & $113$ & $121$ & $129$ & $137$\\
\bottomrule
\end{tabular}\\[6pt]
\end{center}
\caption{Upper bounds for $G(k)$ when $14\le k\le 20$.}\label{tab1}
\end{table}

For comparison, recent work of the second author \cite{Woo2016} delivers the bounds 
$G(14)\le 90$, $G(15)\le 99$, $G(16)\le108$, while rather earlier investigations of Vaughan 
and Wooley \cite{VW2000} obtained $G(17)\le 117$, $G(18)\le 125$, $G(19)\le 134$, 
$G(20)\le 142$. For values of $k$ smaller than $14$, although superior to the bounds of 
\cite{VW2000}, our new methods do not improve on those obtained in \cite{Woo2016}.\par

Two ideas underlie our approach to the theorems above, one old and one new. A novel 
mean value estimate for moments of smooth Weyl sums over sets of minor arcs of 
intermediate and large height is essential for our findings. This new tool is of utility in 
bounding mean values restricted to sets of arcs excluding those of classical major arc type, 
and hence is applicable in pruning problems. A simple but crude version of this idea occurs 
as \cite[Lemma 2.3]{BW2022}, where mean values over sets of major arcs of large height 
are estimated in terms of complete mean values over shortened exponential sums. This 
idea, in turn, has \cite[Lemma 5.6]{LZ2021} as a less flexible and more restricted precursor. 
While a version of \cite[Lemma 2.3]{BW2022} is obtained in Theorem \ref{theorem4.2} 
which applies to lower moments than were accessible hitherto, the treatment of the present 
memoir also delivers analogous bounds for moments restricted to minor arcs. Crucial to our 
applications is the observation that the latter estimates are at their most powerful when the 
associated set of minor arcs is of maximal height relative to the length of the shortened 
exponential sums occurring within our argument. Readers seeking clarity beyond these 
rough and murky remarks would do well to inspect the account in \S5 of the ideas delivering 
Theorem \ref{theorem5.3}.\par

This brings us to the second, much older, idea that we exploit. Minor arc estimates of 
conventional type for smooth Weyl sums over $k$-th powers can be substantially improved 
when their argument lies on an extreme set of minor arcs, rather than on a conventional 
such set. This idea has been utilized previously in work of Heath-Brown \cite{HB1988} 
and Karatsuba \cite{Kar1989}  on fractional parts of $\alp n^k$. A flexible 
analysis using sets of smooth numbers of utility in applications of the circle method can be 
found in \cite{Woo1995}. These improved minor arc estimates can be applied through the 
novel mean value estimates to which we alluded in the previous paragraph, surmounting 
difficulties associated with intermediate sets of arcs that previously obstructed their use. 
The details associated with this plan of attack are described in \S5.\par

We begin the main discourse of this memoir in \S2 by introducing the infrastructure 
required for a discussion of mean values associated with smooth Weyl sums. This section 
already introduces ideas that relate intermediate sets of arcs of differing heights. The 
delicate analysis involved in considering mean values restricted to sets of intermediate arcs 
requires a careful decomposition of smooth Weyl sums, and this we discuss in \S3. Thus 
prepared, we establish our first mean value estimate in \S4, completing the proof of 
Theorem \ref{theorem4.2}. In order to exploit the mean value estimate provided in this 
theorem, we revisit estimates of Weyl type for smooth Weyl sums in \S5, providing in 
Theorem \ref{theorem5.3} an estimate of minor arc type that should be flexible enough for 
future application beyond the present memoir. In \S6 we turn to the application central to 
this paper, namely Waring's problem, and we describe a general analysis. Explicit bounds for 
$G(k)$ are then derived for larger $k$ in \S7, establishing Theorems \ref{theorem1.1} and 
\ref{theorem1.2}. In \S8, we consider intermediate values of $k$ using the tables of 
exponents made available in \cite{VW2000}, and thereby we complete the proof of 
Theorem \ref{theorem1.3}. Finally, in \S9, we briefly outline the consequences of our new 
bounds for problems concerning the representation of almost all positive integers as sums 
of positive integral $k$-th powers.\par

\section{Infrastructure} We initiate the proof of the mean value estimates provided in
Theorems \ref{theorem4.2} and \ref{theorem5.3} by introducing infrastructure necessary 
for the ensuing discussion. A central role is played by the set of $R$-smooth integers not 
exceeding $P$, namely
\[
\mathscr A(P,R)=\{ n\in [1,P]\cap \dbZ:\text{$p|n$ implies $p\le R$}\}.
\]
Here, and throughout this memoir, the letter $p$ is used to denote a prime number. Recall 
the usual convention of writing $e(z)$ for $e^{2\pi iz}$. Then, associated with this set  
$\mathscr A(P,R)$ are the smooth Weyl sum
\[
f(\alp;P,R)=\sum_{x\in \mathscr A(P,R)}e(\alp x^k),
\]
and, for each positive real number $s$, the mean value
\[
U_s(P,R)=\int_0^1|f(\alp;P,R)|^s\d\alp .
\]

\par A real number $\Delta_s$ is referred to as an {\it admissible exponent} (for $k$) if it 
has the property that, whenever $\eps>0$ and $\eta$ is a positive number sufficiently small 
in terms of $\eps$, $k$ and $s$, then whenever $1\le R\le P^\eta$ and $P$ is sufficiently 
large, one has
\[
U_s(P,R)\ll P^{s-k+\Delta_s+\eps}.
\]
Here  and throughout, with $P$ the underlying parameter, the constant implicit in 
Vinogradov's notation may depend on $\eps$, $\eta$, $k$ and $s$. It is easily verified that 
for all positive numbers $s$, one has $\Del_s\ge 0$. It is a simple exercise in interpolation, 
moreover, to confirm that for each $\eta>0$ one has $U_s(P, P^\eta)\gg P^{s/2}$. Thus, 
for all $s>0$ one has
\[\Del_s \ge \max \{ 0,k-s/2\}.
\]
In the opposite direction, one has the trivial upper bound $U_s(P,R) \ll P^s$. Hence 
$\Delta_s=k$ is an admissible exponent. We may therefore suppose that $\Delta_s\le k$, 
and we shall do so whenever this is convenient.\par

We draw a trivial consequence from the definition of an admissible exponent important 
enough that we summarise the conclusion in the form of a lemma.

\begin{lemma}\label{lemma2.1}
Suppose that $\Del_s$ is an admissible exponent for $k$ and that $\eps$ is a positive 
number. Then there exists a positive number $\eta$, depending at most on $\eps$, $k$ and 
$s$, with the following property. Suppose that $P$ is sufficiently large in terms of $\eps$, 
$\eta$, $k$ and $s$, and further that $1\le R\le P^\eta$. Then, uniformly in $1\le Y\le P$, 
one has the bound
\[
U_s(Y,R)\ll P^\eps Y^{s-k+\Del_s}.
\]
\end{lemma}

\begin{proof} Fix $\eps$, $k$ and $s$, so that in our use of Vinogradov's notation we may 
suppress any mention of quantities depending on these numbers, and write 
$\mu_s=s-k+\Del_s$. If we assume that $\Del_s$ is admissible for $k$, there exists a 
positive number $\eta_1$, depending at most on $\eps$, $k$ and $s$, and satisfying 
$\eta_1<\eps$ and the following property. Whenever $X$ is sufficiently large in terms of 
$\eta_1$, say $X\ge X_0(\eta_1)$, and $1\le R\le X^{\eta_1}$, one has 
$U_s(X,R)\ll X^{\mu_s+\eps}$. Now consider a real number $P$ sufficiently large in terms 
of $\eta_1$, and suppose that $1\le Y\le P$. We put $\eta=\eta_1^2/s$ and take $R$ to be 
a real number with $1\le R\le P^\eta$. There are three different regimes for $Y$ that we 
must consider. First, if $Y\le X_0(\eta_1)$, then a trivial estimate yields the bound
\[
U_s(Y,R)\le Y^s\le X_0(\eta_1)^s\ll 1.
\]
Next, when $X_0(\eta_1)<Y\le R^{1/\eta_1}$, the same trivial estimate now reveals that
\[
U_s(Y,R)\le Y^s\le R^{s/\eta_1}\le P^{\eta s/\eta_1}=P^{\eta_1}\le P^\eps.
\]
Finally, when $Y\ge X_0(\eta_1)$ and $R^{1/\eta_1}<Y\le P$, we have $R<Y^{\eta_1}$, 
and then it follows from the above discussion that we have
\[
U_s(Y,R)\ll Y^{\mu_s+\eps}\ll P^\eps Y^{\mu_s}.
\]
By collecting together these estimates, we conclude that the last bound holds uniformly in 
$1\le Y\le P$. This completes the proof of the lemma.
\end{proof}

In order to facilitate concision, from this point onwards we adopt the extended $\eps$, $R$
notation routinely employed by scholars working with smooth Weyl sums while applying the 
Hardy-Littlewood method. Thus, whenever a statement involves the letter $\eps$, then it is 
asserted that the statement holds for any positive real number assigned to $\eps$. Implicit 
constants stemming from Vinogradov or Landau symbols may depend on $\eps$, as well as 
ambient parameters implicitly fixed such as $k$ and $s$. If a statement also involves the 
letter $R$, either implicitly or explicitly, then it is asserted that for any $\eps>0$ there is a 
number $\eta>0$ such that the statement holds uniformly for $2\le R\le P^\eta$. Our 
arguments will involve only a finite number of statements, and consequently we may pass 
to the smallest of the numbers $\eta$ that arise in this way, and then have all estimates in 
force with the same positive number $\eta$. Notice that $\eta$ may be assumed 
sufficiently small in terms of $k$, $s$ and $\eps$.

\par We shall have cause to consider sets of integers, all of whose prime divisors divide a 
fixed integer. In this context, we make use of transparent though disturbing notation, writing 
$u|q^\infty$ to denote that whenever $p$ is a prime and $p|u$, then $p|q$. Then, when 
$q\in \dbN$, we define the set
\[
\mathscr C_q(P,R)=\{ n\in \mathscr A(P,R): n|q^\infty\} ,
\]
consisting of $R$-smooth natural numbers not exceeding $P$ having squarefree kernel 
dividing $q$. We recall that, while $\text{card}(\mathscr A(P,R))\gg_\eta P$ when 
$R\ge P^\eta$, the set $\mathscr C_q(P,R)$ is very thin provided that $q$ is not too large.

\begin{lemma}\label{lemma2.2}
Suppose that $C$ is a positive number. Then, uniformly for positive integers $q$ with 
$q\le P^C$, one has $\text{\rm card}(\mathscr C_q(P,R))\ll P^\eps$.
\end{lemma}

\begin{proof} The desired conclusion is immediate from \cite[Lemma 2.1]{Woo1992}.
\end{proof}

Our interest lies in mean values of $f(\alp,P,R)$ analogous to $U_s(P,R)$, though with 
domains of integration given by intermediate sets of arcs from a Hardy-Littlewood 
dissection. Let $Q$ be a parameter with $1\le Q\le P^{k/2}$. When $q$ is a natural number with $1\le q\le Q$, we define the set of arcs 
$\grM_q(Q,P)$ to be the union of the sets
\[
\grM_{q,a}(Q,P)=\{ \alp\in [0,1): |q\alp-a|\le QP^{-k}\},
\]
with $0\le a\le q$ and $(a,q)=1$, and then put 
\[
\grM(Q,P)=\bigcup_{1\le q\le Q}\grM_q(Q,P).
\]
It is convenient to extend these definitions so that $\grM_q(Q,P)=\emptyset$ when $q>Q$. 
The related dyadically truncated set of arcs $\grN(Q,P)$ may then be defined by
\[
\grN(Q,P)=\grM(Q,P)\setminus \grM(Q/2,P).
\]
Associated with this set are the collections of arcs
\[
\grN_q(Q,P)=\grM_q(Q,P)\setminus \grM_q(Q/2,P).
\]

\par By Dirichlet's approximation theorem, given $\alp\in [0,1)$, there exist $a\in \dbZ$ and 
$q\in \dbN$ with $0\le a\le q\le P^{k/2}$, $(a,q)=1$ and $|q\alp -a|\le P^{-k/2}$. Thus we 
see that $\alp\in \grM(P^{k/2},P)$. Hence, in particular, we have
\[
[0,1)=\bigcup_{j=0}^L\grN(2^{-j}P^{k/2},P),
\]
in which
\begin{equation}
\label{2.1}
L=\left\lfloor \frac{k\log P}{2\log 2}\right\rfloor .
\end{equation}
It therefore follows that
\begin{align*}
U_s(P,R)&=\sum_{j=0}^L\int_{\grN(2^{-j}P^{k/2},P)}|f(\alp;P,R)|^s\d\alp \\
&\ll (\log P)\max_{1\le Q\le P^{k/2}}\int_{\grN(Q,P)}|f(\alp;P,R)|^s\d\alp .
\end{align*}
An important feature of the mean value on the right hand side here is a certain scaling 
property of the associated set $\grN(Q,P)$. We summarise this property in the form of a 
lemma.

\begin{lemma}\label{lemma2.3}
Let $F:\dbR\rightarrow \dbC$ be a $1$-periodic integrable function. Suppose that 
$w\in \dbN$ satisfies the property that $1\le Q\le \tfrac{1}{2}(P/w)^{k/2}$. Then 
whenever $q\in \dbN$ satisfies $(q,w)=1$, one has
\[
\int_{\grM_q(Q,P)}F(\alp w^k)\d\alp =w^{-k}\int_{\grM_q(Q,P/w)}F(\bet)\d\bet .
\]
\end{lemma}

\begin{proof} Let
\[
I=[-q^{-1}QP^{-k},q^{-1}QP^{-k}]\quad \text{and}\quad 
J=[-q^{-1}Qw^kP^{-k},q^{-1}Qw^kP^{-k}].
\]
The hypothesis $Q\le \tfrac12 (P/w)^{k/2}$ ensures that the arcs comprising 
$\grM_q(Q,P/w)$ are disjoint. Since $F$ has period 1, we infer that 
\begin{equation}\label{2.2} 
\int_{\grM_q(Q,P/w)}F(\bet)\d\bet = \sum_{\substack{b=1\\ (b,q)=1}}^q
\int_J F\Big(\frac{b}{q}+\gamma\Big)\d\gamma.
\end{equation}
Likewise, we find that
\begin{align}
\int_{\grM_q(Q,P)}F(\alpha w^k)\d\alpha =&\sum_{\substack{a=1\\ (a,q)=1}}^q\int_I 
F\Big(\Big(\frac{a}{q}+\beta\Big)w^k\Big)\d\beta \notag \\
=&w^{-k}\sum_{\substack{a=1\\ (a,q)=1}}^q
\int_J F\Big(\frac{aw^k}{q}+\gamma\Big)\d\gamma  .\label{2.3}
\end{align}
By hypothesis $(q,w)=1$, whence the mapping $a\mapsto aw^k$ induces a bijection on the 
reduced residue classes modulo $q$. Once again using the hypothesis that $F$ has period 
one, it now follows that the sums on the right hand sides of \eqref{2.2} and \eqref{2.3} 
are equal. This proves the lemma.
\end{proof}

\section{A decomposition of the smooth Weyl sum} We are unable to apply Lemma 
\ref{lemma2.3} directly with $F(\bet)=|f(\bet;P,R)|^s$. However, following a decomposition 
of the smooth Weyl sum $f(\bet;P,R)$, we are able to achieve a conclusion tantamount to 
such an application. Here, the coprimality condition $(q,w)=1$ of Lemma \ref{lemma2.3} 
figures prominently in the analysis. We begin by isolating a part of the smooth Weyl sum 
$f(\alp;P,R)$ in which a large factor $w$ of the argument is available coprime to an auxiliary 
variable $q$. With this objective in mind, we introduce the auxiliary exponential sums
\begin{equation}\label{3.1}
f_q^*(\alp;P,M,R)=\sum_{\substack{v\in \mathscr A(P,R)\\ v>M\\ (v,q)=1}}
\sum_{u\in \mathscr C_q(P/v,R)}e(\alp (uv)^k)
\end{equation}
and
\begin{equation}\label{3.2}
f_q^\dagger(\alp;P,M,R)=\sum_{\substack{v\in \mathscr A(M,R)\\ (v,q)=1}}
\sum_{u\in \mathscr C_q(P/v,R)}e(\alp (uv)^k).
\end{equation}

\begin{lemma}\label{lemma3.1}
Let $q\in\mathbb N$. Then
\[
f(\alp;P,R)=f_q^*(\alp;P,M,R)+f_q^\dagger (\alp;P,M,R).
\]
\end{lemma}

\begin{proof} Consider an integer $x\in \mathscr A(P,R)$, and let $u$ denote the largest 
divisor of $x$ with $u|q^\infty$. Put $v=x/u$. Then either $v\le M$, in which case 
$v\in \mathscr A(M,R)$, or else $v>M$ and $v\in \mathscr A(P,R)$. In both cases, one has 
$x=uv$ with $u\in \mathscr C_q(P/v,R)$ and $(v,q)=1$. The conclusion of the lemma 
follows at once.
\end{proof}

It transpires that the contribution of the exponential sum $f_q^\dagger(\alp;P,M,R)$ is 
easily handled via a trivial estimate.

\begin{lemma}\label{lemma3.2}
Let $Q$ be a parameter with $1\le Q\le P^{k/2}$. Then, whenever $1\le q\le Q$, one has
\[
\int_{\grM_q(Q,P)}|f_q^\dagger(\alp;P,M,R)|^s\d\alp \ll QM^sP^{\eps-k}.
\]
\end{lemma}

\begin{proof} By applying Lemma \ref{lemma2.2} together with a trivial estimate for the 
sum over $v$ in \eqref{3.2}, we see that
\[
|f_q^\dagger (\alp;P,M,R)|\le \sum_{v\le M}\sum_{u\in \mathscr C_q(P/v,R)}1\ll P^\eps M.
\]
Thus, since $\text{mes}(\grM_q(Q,P))\ll QP^{-k}$, we deduce that
\[
\int_{\grM_q(Q,P)}|f_q^\dagger(\alp;P,M,R)|^s\d\alp \ll QP^{-k}(P^\eps M)^s,
\]
and the conclusion of the lemma follows.
\end{proof}

In order to analyse the exponential sum $f_q^*(\alp;P,M,R)$ further, we recall a 
decomposition of the smooth numbers utilised in work  of Vaughan \cite{Vau1989}. In this 
context, we introduce a subset of the smooth numbers $\mathscr A(P,R)$ given by
\[
\mathscr B(M,\pi,R)=\{ v\in \mathscr A(M\pi,R):\text{$v>M$, $\pi|v$ and $\pi'|v$ implies 
$\pi'\ge \pi$}\}.
\] 
Both here and in the remainder of this memoir, we reserve the symbols $\pi$ and $\pi'$ to 
denote prime numbers. We also require the exponential sum
\begin{equation}\label{3.3}
g_{q,\pi}^*(\alp;P,m,R)=\sum_{\substack{w\in \mathscr A(P/m,\pi)\\(w,q)=1}}
\sum_{u\in \mathscr C_q(P/(mw),R)}e(\alp (wu)^k).
\end{equation}

\begin{lemma}\label{lemma3.3} Let $q\in \dbN$. Then whenever $M\ge R$, one has
\[
f_q^*(\alp;P,M,R)=\sum_{\pi\le R}\sum_{\substack{m\in \mathscr B(M,\pi,R)\\ (m,q)=1}}
g_{q,\pi}^*(\alp m^k;P,m,R).
\]
\end{lemma}

\begin{proof} It follows from \cite[Lemma 10.1]{Vau1989} that for each 
$v\in \mathscr A(P,R)$ satisfying $v>M\ge R$, there is a unique triple $(\pi,m,w)$ with 
$v=mw$, $w\in \mathscr A(P/m,\pi)$ and $m\in \mathscr B(M,\pi,R)$. On noting that the 
coprimality conditions $(m,q)=(w,q)=1$ are inherited from the constraint $(v,q)=1$, the 
conclusion of the lemma follows from the definition \eqref{3.1} of $f_q^*(\alp;P,M,R)$. 
\end{proof}

We complete this section by combining the conclusions of Lemmata \ref{lemma3.1}, 
\ref{lemma3.2} and \ref{lemma3.3} so as to obtain a mean value estimate of considerable 
utility. In order to abbreviate notation at this point, we introduce the mean value 
$I_q(M;\grB)$ defined for $\grB$ equal to either $\grM$ or $\grN$ by
\begin{equation}\label{3.4}
I_q(M;\grB)=\sum_{\pi\le R}\sum_{\substack{m\in \mathscr B(M,\pi,R)\\ (m,q)=1}}
\int_{\grB_q(Q,P)}|g_{q,\pi}^*(\alp m^k;P,m,R)|^s\d\alp .
\end{equation}

\begin{lemma}\label{lemma3.4} Let $Q$ be a real number with $1\le Q\le P^{k/2}$, and 
suppose that $s$ is a real number with $s>1$. Then whenever $M\ge R$ and $1\le q\le Q$, 
one has
\[
\int_{\grN_q(Q,P)}|f(\alp;P,R)|^s\d\alp \ll (MR)^{s-1}I_q(M;\grN)+QM^sP^{\eps-k}.
\]
The same conclusion also holds when $\grM$ replaces $\grN$ throughout.
\end{lemma}

\begin{proof} It follows from Lemma \ref{lemma3.1} that when $\alp \in [0,1)$, one has
\[
|f(\alp;P,R)|^s\ll |f_q^*(\alp;P,M,R)|^s+|f_q^\dagger (\alp;P,M,R)|^s.
\]
Moreover, by applying H\"older's inequality in combination with Lemma \ref{lemma3.3}, one 
obtains the bound
\begin{align*}
|f_q^*(\alp;P,M,R)|^s&=\Bigl| \sum_{\pi\le R}
\sum_{\substack{m\in \mathscr B(M,\pi,R)\\ (m,q)=1}}g_{q,\pi}^*(\alp m^k;P,m,R)
\Bigr|^s\\
&\ll (MR)^{s-1}\sum_{\pi\le R}\sum_{\substack{m\in \mathscr B(M,\pi,R)\\ (m,q)=1}}
|g_{q,\pi}^*(\alp m^k;P,m,R)|^s.
\end{align*}
Note that $\grN_q(Q,P)\subseteq \grM_q(Q,P)$. Hence, on integrating over 
$\alp \in \grN_q(Q,P)$ or $\alp \in \grM_q(Q,P)$, the lemma now follows from Lemma 
\ref{lemma3.2}.
\end{proof}

\section{Mean value estimates over intermediate arcs}
The upper bound provided by Lemma \ref{lemma3.4} bounds $f(\alp;P,R)$ in mean, over a 
set of intermediate arcs, in terms of an auxiliary mean value. The latter is susceptible to 
Lemma \ref{lemma2.3}, but the presence of factors in the argument lying in 
$\mathscr C_q(P/(mw),R)$ creates difficulties to which we now attend. In this section, we 
prepare a preliminary mean value using a method that in certain circumstances may be 
enhanced. These enhancements we defer to the next section.\par

We begin with a discussion of the exponential sum $g_{q,\pi}^*(\alp;P,m,R)$. Here, we 
shall find it useful to introduce a modification of the set $\mathscr C_q(P,R)$, namely
\[
\mathscr C_{q,\pi}(P,R)=\{ n\in \mathscr C_q(P,R):\text{$p|n$ implies $p>\pi$}\}.
\]

\begin{lemma}\label{lemma4.1}
One has
\[
g_{q,\pi}^*(\alp;P,m,R)=\sum_{z\in \mathscr C_{q,\pi}(P/m,R)}
\sum_{x\in \mathscr A(P/(mz),\pi)}e(\alp (xz)^k).
\]
\end{lemma}

\begin{proof} On recalling the definition \eqref{3.3} of $g_{q,\pi}^*(\alp;P,m,R)$, we may 
interchange the order of summation to obtain
\[
g_{q,\pi}^*(\alp;P,m,R)=\sum_{u\in \mathscr C_q(P/m,R)}
\sum_{\substack{w\in \mathscr A(P/(mu),\pi)\\ (w,q)=1}}e(\alp (wu)^k).
\]
For each integer $u\in \mathscr C_q(P/m,R)$, there is a unique pair of integers $(y,z)$ 
satisfying $u=yz$, where $y$ has all of its prime divisors no larger than $\pi$, and $z$ has 
no prime divisors less than or equal to $\pi$. Thus, we have $y\in \mathscr C_q(P/m,\pi)$ 
and $z\in \mathscr C_{q,\pi}(P/m,R)$. Making use of this decomposition, we see that
\begin{equation}\label{4.1}
g_{q,\pi}^*(\alp;P,m,R)=\sum_{z\in \mathscr C_{q,\pi}(P/m,R)}
\sum_{y\in \mathscr C_q(P/(mz),\pi)}
\sum_{\substack{w\in \mathscr A(P/(myz),\pi)\\ (w,q)=1}}e(\alp (wyz)^k).
\end{equation}
Notice here that, given any integer $n\in \mathscr A(P/(mz),\pi)$, there are unique integers 
$y$ and $w$ with $n=yw$, and satisfying the condition that $y$ has all of its prime divisors 
amongst those of $q$, and $w$ is coprime with $q$. With such decompositions in mind, we 
recognise that
\[
\sum_{y\in \mathscr C_q(P/(mz),\pi)}
\sum_{\substack{w\in \mathscr A(P/(myz),\pi)\\ (w,q)=1}}e(\gamma (wy)^k)=
\sum_{x\in \mathscr A(P/(mz),\pi)}e(\gamma x^k).
\]
The conclusion of the lemma follows on substituting this relation into \eqref{4.1}.
\end{proof}

We now investigate the mean value $I_q(M;\grB)$ defined in \eqref{3.4} as a prelude to 
the highlight of this section, a mean value estimate for moments of $f(\alp;P,R)$ restricted 
to the set $\grM(Q,P)$. Fix $\grB$ to be either $\grM$ or $\grN$, and fix a real number 
$Q$ with $1\le Q\le \tfrac{1}{2}P^{k/2}R^{-k}$. At this point, we put
\begin{equation}\label{4.2}
M=P(2Q)^{-2/k}R^{-1},
\end{equation}
and we observe that our hypothesis on $Q$ ensures that $M\ge R$. Then, when $\pi\le R$ 
and $m\in \mathscr B(M,\pi, R)$, one has $m\le M\pi\le P(2Q)^{-2/k}$, and thus 
$Q\le \tfrac{1}{2}(P/m)^{k/2}$. The latter condition ensures that the arcs 
$\grM_{q,a}(Q,P/m)$ are disjoint for $0\le a\le q\le Q$ with $(a,q)=1$. Under these 
hypotheses on $Q$ and $m$, therefore, we deduce from \eqref{3.4} via Lemma 
\ref{lemma2.3} that
\begin{equation}\label{4.3}
I_q(M;\grB)=\sum_{\pi\le R}\sum_{\substack{m\in \mathscr B(M,\pi,R)\\ (m,q)=1}}
m^{-k}\int_{\grB_q(Q,P/m)}|g_{q,\pi}^*(\alp;P,m,R)|^s\d\alp .
\end{equation}

\par Observe next that, since $\mathscr C_{q,\pi}(P/m,R)\subseteq \calC_q(P/m,R)$, it 
follows from Lemma \ref{lemma4.1} together with Lemma \ref{lemma2.2} and H\"older's 
inequality that when $s>1$, one has
\begin{align*}
|g_{q,\pi}^*(\alp;P,m,R)|^s&\ll P^\eps \sum_{z\in \mathscr C_{q,\pi}(P/m,R)}
\Bigl| \sum_{x\in \mathscr A(P/(mz),\pi)}e(\alp (xz)^k)\Bigr|^s\\
&\ll P^\eps \sum_{z\in \mathscr A(P/m,R)}|f(\alp z^k;P/(mz),\pi)|^s.
\end{align*}
Write
\[
V_s(\pi,m, z;\grB)=\int_{\grB(Q,P/m)}|f(\alp z^k;P/(mz),\pi)|^s\d\alp .
\]
Then we deduce via \eqref{4.3} that
\begin{align}
\sum_{1\le q\le Q}I_q(M;\grB)&\le \sum_{\pi\le R}\sum_{m\in \mathscr B(M,\pi,R)}m^{-k}
\sum_{1\le q\le Q}\int_{\grB_q(Q,P/m)}|g_{q,\pi}^*(\alp;P,m,R)|^s\d\alp \notag \\
&\ll P^\eps \sum_{\pi\le R}\sum_{m\in \mathscr B(M,\pi,R)}m^{-k}
\sum_{z\in \mathscr A(P/m,R)}V_s(\pi, m,z;\grB).\label{4.4}
\end{align}

The special case of \eqref{4.4} with $\grB=\grM$ combines with Lemma \ref{lemma3.4} to 
deliver the main conclusion of this section. We emphasise that in this statement just as 
elsewhere, we are making use of the extended $\eps$, $R$ convention.

\begin{theorem}\label{theorem4.2}
Suppose that $s$ is a real number with $s\ge 2$ and $\Del_s$ is an admissible exponent. 
Then whenever $Q$ is a real number with $1\le Q\le P^{k/2}$, one has the uniform bound
\[
\int_{\grM(Q,P)}|f(\alp;P,R)|^s\d\alp \ll P^{s-k+\eps}Q^{2\Del_s/k}.
\]
\end{theorem}

\begin{proof} We begin by observing that the conclusion is immediate from the definition of 
an admissible exponent when $\tfrac{1}{2}P^{k/2}R^{-k}<Q\le P^{k/2}$, for in such 
circumstances one has
\[
\int_{\grM(Q,P)}|f(\alp;P,R)|^s\d\alp \le U_s(P,R)\ll P^{s-k+\Del_s+\eps}\ll 
P^{s-k+2\eps}Q^{2\Del_s/k}.
\]
We may therefore suppose henceforth that $1\le Q\le \tfrac{1}{2}P^{k/2}R^{-k}$. In view 
of \eqref{4.2}, one then has also $M\ge R$. For each summand $m$ in the relation 
\eqref{4.4}, one trivially has $\grM(Q,P/m)\subseteq [0,1)$. Thus, by means of a change of 
variable we deduce that
\[
V_s(\pi, m,z;\grM)\le \int_0^1|f(\alp z^k;P/(mz),\pi)|^s\d\alp =U_s(P/(mz),\pi).
\]
We hence infer from Lemma \ref{lemma2.1} and \eqref{4.4} that when 
$s-k+\Del_s\ge 1$, one has
\begin{align*}
\sum_{1\le q\le Q}I_q(M;\grM)&\ll P^\eps \sum_{\pi\le R}
\sum_{m\in \mathscr B(M,\pi,R)}m^{-k}\sum_{z\in \mathscr A(P/m,R)}
\Bigl( \frac{P}{mz}\Bigr)^{s-k+\Del_s}\\
&\ll P^{-k+2\eps }\sum_{\pi\le R}\sum_{m\in \mathscr B(M,\pi,R)}
\Bigl( \frac{P}{m}\Bigr)^{s+\Del_s}\\
&\ll P^{s-k+3\eps}M^{1-s}\Bigl( \frac{P}{M}\Bigr)^{\Del_s}.
\end{align*}

The condition $s-k+\Del_s\ge 1$ is satisfied so long as $s\ge 2$, for as we have already 
observed, it is always the case that $\Del_s\ge k-s/2$. We therefore conclude from Lemma 
\ref{lemma3.4} that
\begin{align*}
\int_{\grM(Q,P)}|f(\alp;P,R)|^s\d\alp &=\sum_{1\le q\le Q}\int_{\grM_q(Q,P)}
|f(\alp;P,R)|^s\d\alp \\
&\ll (MR)^{s-1}\sum_{1\le q\le Q}I_q(M;\grM)+Q^2M^sP^{\eps-k}\\
&\ll P^{s-k+\eps}(P/M)^{\Del_s}+Q^2M^sP^{\eps-k}.
\end{align*}
Thus, on recalling our choice \eqref{4.2} for $M$, we conclude that
\[
\int_{\grM(Q,P)}|f(\alp;P,R)|^s\d\alp \ll P^{s-k+\eps}\left( Q^{2\Del_s/k}
+Q^{2-2s/k}\right). 
\]
The conclusion of the theorem follows on observing that $\Del_s\ge k-s$, whence the first 
term on the right hand side majorises the second.
\end{proof}

We remark that a version of Theorem \ref{theorem4.2} appears as 
\cite[Lemma 2.3]{BW2022}, though in that version the condition $s\ge k+1$ is imposed. 
The proof of that lemma is in many ways more straightforward, with the price being a more 
restrictive constraint on $s$. As we shall see in the next section, the approach that we have 
taken in this memoir also offers the option of retaining minor arc information.

\section{Mean value estimates restricted to minor arcs} The conclusion of Theorem 
\ref{theorem4.2} provides a mean value estimate over an intermediate set of major arcs 
$\grM(Q,P)$. If instead we integrate over the truncated set $\grN(Q,P)$, then we are 
removing the  points from $\grM(Q,P)$ of small height, and the resulting mean value is 
relevant to the estimation of the minor arc contribution. Suppose that 
$1\le Q\le \tfrac{1}{2}P^{k/2}$ and put $X=Q^{2/k}$. Then in very rough terms, one can 
interpret the argument leading to Theorem \ref{theorem4.2} as delivering a bound of the 
flavour
\begin{align*}
\int_{\grM(Q,P)}|f(\alp;P,R)|^s\d\alp &\ll (P/X)^{s-k+\eps}
\int_{\grM(\frac{1}{2}X^{k/2},X)}|f(\alp;X,R)|^s\d\alp \\
&\ll (P/X)^{s-k+\eps}U_s(X,R).
\end{align*}
Our goal now is to obtain an analogous bound of the general shape
\[
\int_{\grN(Q,P)}|f(\alp;P,R)|^s\d\alp \ll (P/X)^{s-k+\eps}
\int_{\grN(\frac{1}{2}X^{k/2},X)}|f(\alp;X,R)|^s\d\alp .
\]
The set $\grN(\frac{1}{2}X^{k/2},X)$ is an extreme set of minor arcs. Here, when 
$\alp$ lies on $\grN(\frac{1}{2}X^{k/2},X)$, it is known that the smooth Weyl sum 
$f(\alp;X,R)$ is $O(X^{1-c/k})$, for a suitable positive number $c$. Since this bound is 
considerably sharper than conventional minor arc bounds for $f(\alp;X,R)$, which would lose 
a factor of roughly $\log k$ in the Weyl exponent, one has rather sharper bounds for
\[
\int_{\grN(Q,P)}|f(\alp;P,R)|^s\d\alp
\]
than were available hitherto, at least when $s$ is fairly large.\par

We begin by deriving a consequence of \cite[Lemma 3.1]{Woo1995}.

\begin{lemma}\label{lemma5.1}
Let $t$ be an even integer, and suppose that the exponent $\Del_t$ is admissible. Then 
whenever $b\in \dbZ$, $r\in \dbN$ and $(b,r)=1$, one has
\[
f(\alp;P,R)\ll r^\eps P^{1+\eps}\left( P^{\Del_t}\left( \Tet^{-1}+P^{-k/2}+P^{-k}\Tet 
\right) \right)^{2/t^2} +P^{1/2+\eps},
\]
in which we write $\Tet=r+P^k|r\alp-b|$.
\end{lemma}

\begin{proof} Suppose that $\tfrac{1}{2}<\lam <1$, $M=P^\lam$ and $\alp\in \dbR$. 
Suppose further that $a\in \dbZ$ and $q\in \dbN$ satisfy $(a,q)=1$ and 
$|\alp-a/q|\le 1/q^2$. Then \cite[Lemma 3.1]{Woo1995} establishes that for all even 
natural numbers $t$ and $w$, one has
\begin{align*}
f(\alp;P,R)\ll &\, q^\eps P^{1+\eps} \left( M^{\Del_w}(P/M)^{\Del_t}
\left( q^{-1}+M^{-k}+(P/M)^{-k}+qP^{-k}\right) \right)^{2/(tw)}\\
&\,+M,
\end{align*}
We take $w=t$ and $\lam=\tfrac{1}{2}+\del$, for a small fixed positive number $\del$. 
Thus
\[
f(\alp;P,R)\ll q^\eps P^{1+k\del} \left( P^{\Del_t}\left( q^{-1}+P^{-k/2}+qP^{-k}\right) 
\right)^{2/t^2}+P^{1/2+\del}.
\]
We now apply a standard transference principle (see \cite[Lemma 14.1]{Woo2015}) to see 
that the same conclusion holds for all $b\in \dbZ$ and $r\in \dbN$ with $(b,r)=1$ when we 
replace $q$ by $\Tet=r+P^k|r\alp-b|$ throughout. The conclusion of the lemma therefore 
follows, since $\del$ may be taken arbitrarily small.
\end{proof}

The most powerful consequences of Lemma \ref{lemma5.1} are made available by applying 
Dirichlet's approximation theorem to obtain integers $b$ and $r$ with $(b,r)=1$ and 
$1\le r\le P^{k/2}$ for which $|r\alp-b|\le P^{-k/2}$. In such circumstances, Lemma 
\ref{lemma5.1} is most effective when $\alp$ satisfies the condition that $r>cP^{k/2}$, for 
some fixed  $c>0$. One then has $f(\alp;P,R)\ll P^{1-\tau(t,k)+\eps}+P^{1/2+\eps}$, 
where
\[
\tau(t,k)=\frac{k-2\Del_t}{t^2}.
\]
Since $\Del_t\ge \max \{ k-t/2,0\}$, one sees that
\[
\tau(t,k)\le \min \left\{ \frac{t-k}{t^2},\frac{k}{t^2}\right\}\le \frac{1}{4k},
\]
and thus our estimate for $f(\alpha;P,R)$ simplifies to $f(\alp;P,R)\ll P^{1-\tau(t,k)+\eps}$. 
To extract the most from this bound, we introduce the number 
\begin{equation}\label{5.2}
\tau(k)=\max_{w\in \dbN}\frac{k-2\Del_{2w}}{4w^2},
\end{equation}
and then have
\begin{equation}\label{5.1}
f(\alp;P,R)\ll P^{1-\tau(k)+\eps}.
\end{equation}
The number $\tau(k)$ will be of significance in the argument below. It appears also in 
slightly different guises in work of Karatsuba \cite{Kar1989} and Heath-Brown 
\cite{HB1988}.

\medskip

We now return to the rescaling argument underlying the work of \S4. In this context, we 
introduce an auxiliary exponent. Suppose that $s$ is a real number with $s\ge 2$, and that 
the exponents $\Del_u$ are admissible for $2\le u\le s$. We define
\begin{equation}\label{5.3}
\Del_s^*=\min_{0\le t\le s-2}\left( \Del_{s-t}-t\tau(k)\right) ,
\end{equation}
and refer to $\Del_s^*$ as an {\it admissible exponent for minor arcs}.

\begin{theorem}\label{theorem5.2}
Suppose that $s\ge 2$, and that $\Del_s^*$ is an admissible exponent for minor arcs. Then 
whenever $1\le Q\le \tfrac{1}{2}P^{k/2}R^{-k}$, one has the uniform bound
\[
\int_{\grN(Q,P)}|f(\alp;P,R)|^s\d\alp \ll P^{s-k+\eps}Q^{2\Del_s^*/k}.
\]
\end{theorem}

\begin{proof} We again fix $M$ according to equation \eqref{4.2}, and we recall from 
\eqref{4.4} that when $1\le Q\le \tfrac{1}{2}P^{k/2}R^{-k}$, one has
\begin{equation}\label{5.4}
\sum_{1\le q\le Q}I_q(M;\grN)\ll P^\eps \sum_{\pi\le R}
\sum_{m\in \mathscr B(M,\pi,R)}m^{-k}\sum_{z\in \mathscr A(P/m,R)}V_s(\pi,m,z;\grN),
\end{equation}
where
\begin{equation}\label{5.5}
V_s(\pi,m,z;\grN)=\int_{\grN(Q,P/m)}|f(\alp z^k;P/(mz),\pi)|^s\d\alp .
\end{equation}

\par We apply Lemma \ref{lemma5.1} to estimate $f(\alp z^k;P/(mz),\pi)$ when 
$\alp \in \grN(Q,P/m)$. In the latter circumstances, one has $\alp\in \grM(Q,P/m)\setminus 
\grM(Q/2,P/m)$. Thus, there exist integers $a$ and $q$ with $0\le a\le q\le Q$ and 
$(a,q)=1$ for which one has $|q\alp -a|\le Q(P/m)^{-k}$, and either $q>Q/2$ or 
$|q\alp-a|>\tfrac{1}{2}Q(P/m)^{-k}$. Consider a fixed integer $z\in \mathscr A(P/m,R)$. 
Then as a consequence of these relations, if we put
\[
r=\frac{q}{(q,z^k)}\quad \text{and}\quad b=\frac{az^k}{(q,z^k)},
\]
then we find that $(b,r)=1$ with $r\le Q$ and $|r(\alp z^k)-b|\le Q(P/(mz))^{-k}$. 
Moreover, one has either $r>\tfrac{1}{2}Qz^{-k}$ or 
$|r(\alp z^k)-b|>\tfrac{1}{2}Q(P/m)^{-k}$. Thus, in particular,
\[
\tfrac{1}{2}Qz^{-k}<r+\Bigl(\frac{P}{mz}\Bigr)^k|r(\alp z^k)-b|\le 2Q.
\]
We therefore deduce from Lemma \ref{lemma5.1} that whenever $t$ is an even integer, 
then
\begin{align*}
f(\alp z^k;P/(mz),\pi)\ll &\,Q^\eps \Bigl( \frac{P}{mz}\Bigr)^{1+\eps} \left( 
\Bigl( \frac{P}{mz}\Bigr)^{\Del_t}\Bigl( \frac{z^k}{Q}+\Bigl( \frac{mz}{P}\Bigr)^{k/2}+
Q\Bigl( \frac{mz}{P}\Bigr)^k \Bigr) \right)^{2/t^2}\\
&\, +\Bigl( \frac{P}{mz}\Bigr)^{1/2+\eps}.
\end{align*}
We choose $t=2w$ to correspond to the maximum in the definition of $\tau=\tau(k)$ in 
\eqref{5.2}, and recall from \eqref{4.2} that $Q=\tfrac{1}{2}(P/(MR))^{k/2}$. Then, when 
$M<m\le MR$ and $\alp\in \grN(Q,P/m)$, we conclude that
\[
f(\alp z^k;P/(mz),\pi)\ll \Bigl( \frac{P}{mz}\Bigr)^{1/2+\eps}+P^\eps 
\Bigl( \frac{P}{m}\Bigr)^{1-\tau}z^{-1+2(k-\Del_t)/t^2}.
\]
Since $\Del_t\ge k-t/2$ and $t\ge 2$, we arrive at the upper bound
\begin{equation}\label{5.6}
\sup_{\alp\in \grN(Q,P/m)}|f(\alp z^k;P/(mz),\pi)|\ll P^\eps (P/m)^{1-\tau}z^{-1/2}.
\end{equation}

\par We now return to the mean value $V_s(\pi,m,z;\grN)$ defined in \eqref{5.5}. Let $t$ 
and $v$ be non-negative integers with $s=t+v$. Then it follows from \eqref{5.6} that
\[
V_s(\pi,m,z;\grN)\ll P^\eps (P/m)^{t(1-\tau)}\int_0^1|f(\alp z^k;P/(mz),\pi)|^v\d\alp .
\]
A change of variable therefore combines with Lemma \ref{lemma2.1} to show that
\begin{align*}
V_s(\pi,m,z;\grN)&\ll P^\eps (P/m)^{t(1-\tau)}U_v(P/(mz),\pi)\\
&\ll P^{2\eps}(P/m)^{t(1-\tau)}(P/(mz))^{v-k+\Del_v}.
\end{align*}
Since $\Del_v\ge k-v/2$, we see that $v-k+\Del_v\ge 1$ whenever $v\ge 2$. On recalling 
the definition \eqref{5.3} of $\Del_s^*$, therefore, and noting that $v=s-t$, we discern 
that
\[
V_s(\pi,m,z;\grN)\ll z^{-1}P^\eps (P/m)^{s-k+\Del_s^*}.
\]
On substituting this upper bound into \eqref{5.4}, we find that
\begin{align}
\sum_{1\le q\le Q}I_q(M;\grN)&\ll P^\eps \sum_{\pi\le R}\sum_{M<m\le MR}m^{-k}
(P/m)^{s-k+\Del_s^*}\sum_{1\le z\le P/m}z^{-1}\notag \\
&\ll P^{s-k+2\eps}M^{1-s}(P/M)^{\Del_s^*}.\label{5.7}
\end{align}

\par We next appeal to Lemma \ref{lemma3.4}, proceeding just as in the conclusion of the 
proof of Theorem \ref{theorem4.2}. Thus, making use of the bound \eqref{5.7}, we obtain
\begin{align*}
\int_{\grN(Q,P)}|f(\alp;P,R)|^s\d\alp &=\sum_{1\le q\le Q}\int_{\grN_q(Q,P)}
|f(\alp;P,R)|^s\d\alp \\
&\ll (MR)^{s-1}\sum_{1\le q\le Q}I_q(M;\grN)+Q^2M^sP^{\eps-k}\\
&\ll P^{s-k+\eps}(P/M)^{\Del_s^*}+Q^2M^sP^{\eps-k}.
\end{align*}
Hence, on recalling the choice \eqref{4.2} for $M$, we conclude that
\begin{equation}\label{5.8}
\int_{\grN(Q,P)}|f(\alp;P,R)|^s\d\alp \ll P^{s-k+\eps}(Q^{2-2s/k}+Q^{2\Del_s^*/k}).
\end{equation}
We have observed already that $\tau(k)\le 1/(4k)$. Thus, since $\Del_{s-t}\ge k-(s-t)$, 
one sees that for some integer $t$ satisfying $0\le t\le s-2$ (the integer $t$ associated with 
the definition \eqref{5.3} of $\Del_s^*$), one has
\[
\frac{2}{k}\Del_s^*\ge \frac{2}{k}\Bigl( k-(s-t)-\frac{t}{4k}\Bigr) \ge 2-\frac{2s}{k}.
\]
The desired conclusion is therefore immediate from \eqref{5.8}.
\end{proof}

This theorem may be exploited to obtain a bound for minor arc contributions of considerable 
utility in applications of the circle method. In this context, we introduce the set of minor arcs 
$\grm(Q)=\grm(Q,P)$ given by $\grm(Q)=[0,1]\setminus \grM(Q,P)$. We also abbreviate 
the major arcs $\grM(Q,P)$ simply to $\grM(Q)$ in circumstances where the implicit second 
parameter is equal to $P$ and brevity is to be prized above full disclosure. 

\begin{theorem}\label{theorem5.3}
Let $s\ge 2$ and suppose that $\Del_s^*$ is an admissible exponent for minor arcs 
satisfying $\Del_s^*<0$. Let $\tet$ be a positive number with $\tet\le k/2$. Then whenever 
$P^\tet\le Q\le P^{k/2}$, one has the bound
\[
\int_{\grm(Q)}|f(\alp;P,R)|^s\d\alp \ll_\tet P^{s-k}Q^{\eps-2|\Del_s^*|/k}.
\]
\end{theorem}

\begin{proof} Write
\[
J=\left\lceil \frac{\log (P^{k/2}/Q)}{\log 2}\right\rceil \quad \text{and}\quad 
J_0=\left\lceil \frac{\log (2R^k)}{\log 2}\right\rceil .
\]
We begin by observing that, since $\grm(Q)=[0,1]\setminus \grM(Q,P)$, we have
\[
\grm(Q)\subseteq \bigcup_{j=0}^J \grN(2^{-j}P^{k/2},P).
\]
When $J_0<j\le J$, it follows from Theorem \ref{theorem5.2} that
\begin{align}
\int_{\grN(2^{-j}P^{k/2},P)}|f(\alp;P,R)|^s\d\alp &\ll P^{s-k+\eps}\left( 2^{-j}P^{k/2}
\right)^{2\Del_s^*/k}\notag \\
&\ll P^{s-k+\eps}Q^{-2|\Del_s^*|/k}.\label{5.9}
\end{align}

\par Meanwile, when $0\le j\le J_0$, we may apply the argument underlying the proof of 
Theorem \ref{theorem5.2}. Thus, when $\alp\in \grN(2^{-j}P^{k/2},P)$, there exist 
$b\in \dbZ$ and $r\in \dbN$ with $(b,r)=1$, $r\le 2^{-j}P^{k/2}$ and 
$|r\alp-b|\le 2^{-j}P^{-k/2}$. Since $\alp\not\in \grM(2^{-j-1}P^{k/2},P)$, we have
\[
P^{k/2}R^{-k}\ll 2^{-1-j}P^{k/2}\le r+P^k|r\alp-b|\ll P^{k/2}.
\]
By Lemma \ref{lemma5.1} and \eqref{5.1}, we now have 
$f(\alp;P,R)\ll P^{1-\tau(k)+\eps}$. With $s=t+v$, and $t$ and $v$ defined as in the proof 
of Theorem \ref{theorem5.2}, we therefore infer that
\begin{align*}
\int_{\grN(2^{-j}P^{k/2},P)}|f(\alp;P,R)|^s\d\alp &\ll (P^{1-\tau(k)+\eps})^t
\int_0^1|f(\alp;P,R)|^v\d\alp \\
&\ll P^{s-k+\eps}P^{\Del_v-t\tau(k)} \\
&\ll P^{s-k+\eps}Q^{-2|\Del_s^*|/k}.
\end{align*}
On combining this estimate with \eqref{5.9}, we see that
\begin{align*}
\int_{\grm(Q)}|f(\alp;P,R)|^s\d\alp &\ll \sum_{j=0}^J\int_{\grN(2^{-j}P^{k/2},P)}
|f(\alp;P,R)|^s\d\alp \\
&\ll P^{s-k+\eps}Q^{-2|\Del_s^*|/k}.
\end{align*}
Since $Q\ge P^\tet$ and $\tet>0$, it suffices to recall the conventions concerning the use 
of $\varepsilon$ and $\eta$ to complete the proof of the theorem.
\end{proof}

\section{The treatment of $G(k)$ in general terms}
Our proofs of Theorems \ref{theorem1.1} and \ref{theorem1.2} are largely routine given 
the flexible nature of Theorem \ref{theorem5.3}, so we may be concise in our exposition. 
We begin with a pruning argument that extends the range of $Q$ in Theorem 
\ref{theorem5.3} from a power of $P$ to an arbitrarily slowly growing function of $P$.

\begin{theorem}\label{theorem6.1}
Suppose that $k\ge 3$, $s\ge 2k+3$ and $\Del_s^*$ is an admissible exponent for minor 
arcs with $\Del_s^*<0$. Let $\nu$ be any positive number with 
\[
\nu<\min \left\{ \frac{2|\Del_s^*|}{k}, \frac{1}{6k}\right\} .
\]
Then, when $1\le Q\le P^{k/2}$, one has the uniform bound
\[
\int_{\grm(Q)}|f(\alp;P,R)|^s\d\alp \ll P^{s-k}Q^{-\nu}.
\]
\end{theorem}

\begin{proof} In view of the conclusion of Theorem \ref{theorem5.3}, it suffices to 
consider values of $Q$ with $1\le Q\le P^\tet$, where $\tet$ is a fixed positive number 
small in terms of $k$ and $s$. We assume in particular that $\tet<1/k$, whence for 
$k\ge 3$ one has
\begin{equation}\label{6.1}
\frac{3}{4}+\frac{\tet}{8}<1-\frac{1}{2k}.
\end{equation}
Our starting point is the observation that, as a consequence of Theorem \ref{theorem5.3},
\begin{align}
\int_{\grm(Q)}|f(\alp;P,R)|^s\d\alp &=\int_{\grm(P^\tet)}|f(\alp;P,R)|^s\d\alp +
\int_{\grm(Q)\setminus \grm(P^\tet)}|f(\alp;P,R)|^s\d\alp \notag \\
&\ll P^{s-k}(P^\tet)^{\eps-2|\Del_s^*|/k}+\int_{\grm(Q)\setminus \grm(P^\tet)}
|f(\alp;P,R)|^s\d\alp \notag \\
&\ll P^{s-k}Q^{-\nu}+\int_{\grM(P^\tet)\setminus \grM(Q)}|f(\alp;P,R)|^s\d\alp .
\label{6.2}
\end{align} 

\par When $a\in \dbZ$ and $q\in \dbN$ satisfy $0\le a\le q\le \tfrac{1}{2}P^{k/2}$ and 
$(a,q)=1$, the intervals $\grM_{q,a}(\tfrac{1}{2}P^{k/2},P)$ are disjoint, and for 
$\alp\in \grM_{q,a}(\tfrac{1}{2}P^{k/2},P)$ we put
\[
\Ups(\alp)=(q+P^k|q\alp-a|)^{-1}.
\]
Meanwhile, for $\alp\in [0,1)\setminus \grM(\tfrac{1}{2}P^{k/2},P)$ we put 
$\Ups(\alp)=0$. This defines a function $\Ups: [0,1)\rightarrow [0,1]$. By 
\cite[Lemma 7.2]{VW1991}, we find that when
\[
2\le R\le M\le P,\quad |q\alp-a|\le M/(k(2P)^kR)\quad \text{and}\quad (a,q)=1,
\]
one has
\[
f(\alp;P,R)\ll q^\eps L^3\left( P\Ups(\alp)^{1/(2k)}+(PMR)^{1/2}+q^{1/4}P(R/M)^{1/2}
\right) 
\]
where $L$ is defined by \eqref{2.1}. But on taking $M=P^{(2+\tet)/4}$ and recalling 
\eqref{6.1}, we see that when $q\le P^\tet$ one has
\[
q^\eps L^3\left( (PMR)^{1/2}+q^{1/4}P(R/M)^{1/2}\right) \ll 
P^{\frac{3}{4}+\frac{\tet}{8}+\eps}R\ll P^{1-1/(2k)}.
\]
It follows that whenever $\alp\in \grM(P^\tet,P)\setminus \grM (Q,P)$, one has the bound
\begin{equation}\label{6.3}
f(\alp;P,R)\ll PL^3\Ups(\alp)^{-\eps+1/(2k)}+P^{1-\tau(k)+\eps}.
\end{equation}

\par We now put $s=t+v$, where $t$ and $v$ are chosen in accordance with the definition 
\eqref{5.3} of $\Del_s^*$, just as in the proof of Theorem \ref{theorem5.2}. Thus, by 
substituting \eqref{6.3} into \eqref{6.2}, we obtain the bound
\begin{equation}\label{6.4}
\int_{\grm(Q)}|f(\alp;P,R)|^s\d\alp \ll P^{s-k}Q^{-\nu} +P^\eps T_1+(PL^3)^tT_2,
\end{equation}
where
\begin{equation}\label{6.5}
T_1=\left(P^{1-\tau(k)}\right)^t\int_0^1|f(\alp;P,R)|^v\d\alp 
\end{equation}
and
\begin{equation}\label{6.6}
T_2=\int_{\grM (P^\tet)\setminus \grM(Q)}\Ups(\alp)^{-\eps+t/(2k)}
|f(\alp;P,R)|^v\d\alp .
\end{equation}

\par As in the proof of Theorem \ref{theorem5.3}, it is apparent from \eqref{6.5} that
\[
T_1\ll \left(P^{1-\tau(k)}\right)^tP^{v-k+\Del_v+\eps}\ll P^{s-k-|\Del_s^*|+\eps}.
\]
Thus we obtain the estimate
\begin{equation}\label{6.7}
P^\eps T_1\ll P^{s-k}Q^{-\nu}.
\end{equation}

\par Meanwhile, an application of H\"older's inequality to \eqref{6.6} reveals that
\begin{equation}\label{6.8}
T_2\le T_3^{(v-2)/(s-2)}T_4^{t/(s-2)},
\end{equation}
where
\[
T_3=\int_{\grm(Q)}|f(\alp;P,R)|^s\d\alp 
\]
and
\begin{equation}\label{6.9}
T_4=\int_{\grM(P^\tet)\setminus \grM(Q)}\Ups(\alp)^{-\eps+(s-2)/(2k)}
|f(\alp;P,R)|^2\d\alp .
\end{equation}
On substituting \eqref{6.7} and \eqref{6.8} into \eqref{6.4}, we obtain the estimate
\[
T_3\ll P^{s-k}Q^{-\nu}+Q^\eps (PL^3)^tT_3^{1-t/(s-2)}T_4^{t/(s-2)},
\]
whence
\begin{equation}\label{6.10}
\int_{\grm(Q)}|f(\alp;P,R)|^s\d\alp \ll P^{s-k}Q^{-\nu}+(Q^\eps PL^3)^{s-2}T_4.
\end{equation}
Thus it remains only to bound the mean value $T_4$.\par

When $\alp \in \dbR$, it follows from Dirichlet's approximation theorem that there exist 
$a\in \dbZ$ and $q\in \dbN$ with $(a,q)=1$, $q\le Q^{-1}P^k$ and 
$|q\alp-a|\le QP^{-k}$. When $\alp\in \grM(P^\tet)\setminus \grM(Q)$, moreover, one has
$q+P^k|q\alp -a|>Q$, and hence $\Ups(\alp)<Q^{-1}$. We therefore deduce from 
\eqref{6.9} that when $s\ge 2k+3$, we have the bound
\[
T_4\ll Q^{\eps -1/(4k)}\int_{\grM(P^\theta)}\Ups(\alp)^{1+1/(4k)}|f(\alp;P,R)|^2\d\alp .
\]
The mean value on the right hand side here is amenable to \cite[Lemma 11.1]{PW2014}, a 
pruning lemma that refines earlier work of the first author \cite[Lemma 2]{Bru1988}. Thus, 
we obtain the estimate $T_4\ll Q^{\eps -1/(4k)}P^{2-k}$. After substituting this bound into 
\eqref{6.10}, we infer that
\[
\int_{\grm(Q)}|f(\alp;P,R)|^s\d\alp \ll P^{s-k}Q^{-\nu}+Q^{\eps -1/(4k)}L^{3s}P^{s-k}.
\]
The desired conclusion therefore follows provided that $Q>L^{60ks}$, since then
\[
L^{3s}Q^{\eps-1/(4k)}\le Q^{\eps -1/(5k)}\left( L^{60ks}Q^{-1}\right)^{1/(20k)}<
Q^{-\nu}.
\]

\par At this point we are reduced to the scenario in which one has $Q\le L^{60ks}$. In this 
range for $Q$, we appeal to \cite[Lemma 8.5]{VW1991}. Let $A>0$ be fixed. Then the 
latter lemma shows that when $a\in \dbZ$ and $q\in \dbN$ satisfy $(a,q)=1$ and 
$q\le L^A$, one has the upper bound
\[
f(\alp;P,R)\ll P\Ups(\alp)^{-\eps+1/k}+
P\exp \left( -c(\log P)^{1/2}\right) (1+P^k|\alp-a/q|),
\]
in which $c=c(A)>0$. When $\alp\in \grM(L^{60ks})\setminus \grM(Q)$, one has
\[
Q<q+P^k|q\alp-a|\le 2L^{60ks}.
\]
In such circumstances, therefore, we have
\[
f(\alp;P,R)\ll P\Ups(\alp)^{-\eps+1/k}+PL^{-60ks}\ll P\Ups(\alp)^{1/(2k)}Q^{-1/(3k)}.
\]
Write
\[
T_5=\int_{\grM(L^{60ks})\setminus \grM(Q)}|f(\alp;P,R)|^s\d\alp .
\]
Then we deduce that when $s\ge 2k+3$, one has
\begin{align*}
T_5&\ll (PQ^{-1/(3k)})^{s-2}\int_{\grM(L^{60ks})}\Ups(\alp)^{(s-2)/(2k)}
|f(\alp;P,R)|^2\d\alp \\
&\ll P^{s-2}Q^{-1/2}\int_{\grM(L^{60ks})}\Ups(\alp)^{1+1/(2k)}|f(\alp;P,R)|^2\d \alp.
\end{align*}
Observe that $\grm(Q)\setminus \grm(L^{60ks})=\grM(L^{60ks})\setminus \grM(Q)$. 
Then, again employing \cite[Lemma 11.1]{PW2014}, we conclude that
\[
\int_{\grm(Q)\setminus \grm(L^{60ks})}|f(\alp;P,R)|^s\d\alp =T_5\ll P^{s-k}Q^{-1/2}.
\]
Hence, on applying the conclusion of the theorem already established when 
$Q\ge L^{60ks}$, we obtain
\begin{align*}
\int_{\grm(Q)}|f(\alp;P,R)|^s\d\alp &=\int_{\grm(L^{60ks})}|f(\alp;P,R)|^s\d\alp +T_5\\
&\ll P^{s-k}(L^{60ks})^{-\nu}+P^{s-k}Q^{-1/2}\\
&\ll P^{s-k}Q^{-\nu}.
\end{align*}
The conclusion of the theorem therefore follows also in this last case with 
$1\le Q\le L^{60ks}$, and thus the proof of the theorem is complete.
\end{proof}

We are now equipped to bound the quantity $G(k)$ relevant to Waring's problem. We 
assume that we have available an admissible exponent $\Del_u$ for each positive number 
$u$. Then, when $k\ge 4$, we define $\tau(k)$ as in \eqref{5.2},
and we also put
\begin{equation}\label{6.11}
G_0(k)=\min_{v\ge 2}\left( v+\frac{\Del_v}{\tau(k)}\right) .
\end{equation}
Also, when $s\in \dbN$, we write $R_{s,k}(n)$ for the number of solutions of the equation 
\begin{equation}\label{6.12}
x_1^k+\ldots +x_s^k=n,
\end{equation}
with $x_i\in \dbN$.

\begin{theorem}\label{theorem6.2}
Suppose that $k\ge 4$ and $s\ge \max \{ \lfloor G_0(k)\rfloor +1,2k+3\}$. Then provided 
that the integer $n$ is sufficiently large in terms of $k$ and $s$, and for each natural 
number $q$ the congruence
\[
x_1^k+\ldots +x_s^k\equiv n\mmod{q}
\]
possesses a solution with $(x_1,q)=1$, one has $R_{s,k}(n)\gg n^{s/k-1}$. In particular, 
when $k$ is not a power of $2$ one has $G(k)\le \max \{ \lfloor G_0(k)\rfloor +1,2k+3\}$, 
and when $k$ is a power of $2$ one has instead 
$G(k)\le \max \{ \lfloor G_0(k)\rfloor +1,4k\}$.
\end{theorem}

\begin{proof} We first address the claimed asymptotic lower bound 
$R_{s,k}(n)\gg n^{s/k-1}$, the final conclusions of the theorem following from the 
standard theory associated with local solubility in Waring's problem (see 
\cite[Theorem 4.6]{Vau1997}, for example). Consider a natural number $n$ sufficiently 
large in terms of $k$ and $s$. Let $P=n^{1/k}$ and $R=P^\eta$, where $\eta>0$ is 
sufficiently small, in a manner to be specified in due course. We denote by $r_{s,k}(n)$ the 
number of representations of $n$ in the form \eqref{6.12} with $x_i\in \mathscr A(P,R)$ 
$(1\le i\le s)$, so that $R_{s,k}(n)\ge r_{s,k}(n)$. By orthogonality, one has
\[
r_{s,k}(n)=\int_0^1f(\alp;P,R)^se(-n\alp)\d\alp .
\]
We put $Q=L^{1/15}$, and we specify $\eta$ to be sufficiently small in the context of the 
(finitely many) admissible exponents that must be discussed in determining $\tau(k)$ and 
$G_0(k)$. We make use of a simplified Hardy-Littlewood dissection. Thus, we take $\grK$ 
to be the union of the arcs
\[
\grK(q,a)=\{ \alp \in [0,1):|\alp-a/q|\le QP^{-k}\},
\]
with $0\le a\le q\le Q$ and $(a,q)=1$, and then put $\grk=[0,1)\setminus \grK$. Thus, by 
the triangle inequality, we have
\begin{equation}\label{6.13}
r_{s,k}(n)=\int_\grK f(\alp;P,R)^se(-n\alp)\d\alp +O\biggl( \int_\grk |f(\alp;P,R)|^s\d\alp 
\biggr) .
\end{equation} 

We first handle the contribution of the minor arcs $\grk$ within \eqref{6.13}. Suppose that 
$s\ge \max \{ \lfloor G_0(k)\rfloor +1,2k+3\}$, and recall \eqref{5.3} and \eqref{6.11}. 
Then there exists a positive number $v$ with $v\ge 2$ and an admissible exponent 
$\Del_v$ for which the exponent $\Del_s^*$ is admissible for minor arcs, where
\[
\Del_s^*=\Del_v-(s-v)\tau(k)=-\tau(k) \left( s-G_0(k)\right)<0.
\]  
Put $\nu=\min \{ |\Del_s^*|/k, 1/(18k)\}$. Then we see from Theorem \ref{theorem6.1} 
that
\[
\int_{\grm (Q)}|f(\alp;P,R)|^s\d\alp \ll P^{s-k}Q^{-\nu}=P^{s-k}L^{-\nu/15}.
\]
Finally, since $\grk\subseteq \grm(Q)$, we may conclude thus far that
\begin{equation}\label{6.14}
\int_\grk |f(\alp;P,R)|^s\d\alp \le \int_{\grm (Q)}|f(\alp;P,R)|^s\d\alp \ll 
P^{s-k}L^{-\nu/15}.
\end{equation} 

\par Next we attend to the contribution of the major arcs $\grK$. Suppose that 
$\alp\in \grK(q,a)\subseteq \grK$. The standard theory of smooth Weyl sums (see 
\cite[Lemma 5.4]{Vau1989}) shows that there is a positive number $c=c(\eta)$ such that
\[
f(\alp;P,R)=cq^{-1}S(q,a)v(\alp-a/q)+O(PL^{-1/4}),
\]
wherein
\[
S(q,a)=\sum_{r=1}^qe(ar^k/q)\quad \text{and}\quad v(\bet)=\frac{1}{k}\sum_{m\le n}
m^{-1+1/k}e(\bet m).
\]
Since $\grK$ has measure $O(Q^3n^{-1})$, we see that
\begin{equation}\label{6.15}
\int_\grK f(\alp;P,R)^se(-n\alp)\d\alp =c^s\grJ(n,Q)\grS(n,Q)+O(P^{s-k}Q^3L^{-1/4}),
\end{equation}
where
\[
\grJ(n,X)=\int_{-X/n}^{X/n}v(\bet)^se(-\bet n)\d\bet 
\]
and
\[
\grS(n,X)=\sum_{1\le q\le X}\sum^q_{\substack{a=1\\ (a,q)=1}}q^{-s}S(q,a)^se(-na/q).
\]
Notice that since $Q=L^{1/15}$, the error term in \eqref{6.15} is $O(P^{s-k}L^{-1/20})$. 
Familiar estimates from the theory of Waring's problem (see 
\cite[Chapters 2 and 4]{Vau1997}) show that under the hypotheses on $s$ at hand,
\[
\grS(n,X)=\grS(n)+O(X^{-1/k}),
\]
where
\[
\grS(n)=\sum_{q=1}^\infty \sum^q_{\substack{a=1\\ (a,q)=1}}q^{-s}S(q,a)^se(-na/q).
\]
Thus, in particular, subject to the hypotheses of the statement of the theorem, one has 
$\grS(n,X)\gg 1$. Likewise, one finds that
\[
\grJ(n,X)=\frac{\Gam(1+1/k)^s}{\Gam(s/k)}n^{s/k-1}+O(n^{s/k-1}X^{-1/k}).
\]
Hence, again under the hypotheses of the statement of the theorem, we deduce from 
\eqref{6.15} that
\begin{equation}\label{6.16}
\int_\grK f(\alp;P,R)^se(-n\alp)\d\alp =c^s\grS(n)
\frac{\Gam(1+1/k)^s}{\Gam(s/k)}n^{s/k-1}+o(n^{s/k-1}).
\end{equation}

\par On substituting \eqref{6.14} and \eqref{6.16} into \eqref{6.13}, we conclude that
\[
r_{s,k}(n)=c^s\grS(n)\frac{\Gam(1+1/k)^s}{\Gam(s/k)}n^{s/k-1}+o(n^{s/k-1}),
\]
whence $R_{s,k}(n)\ge r_{s,k}(n)\gg n^{s/k-1}$. This completes the proof of the 
asymptotic lower bound asserted in the statement of the theorem, subject of course to 
the associated hypotheses on $s$, and, when $s<4k$ and $k$ is a power of $2$, the 
hypothesis on local solubility. Since we have already confirmed the remaining assertions of 
the theorem, subject to validity of this asymptotic lower bound, the proof of the theorem is 
complete.
\end{proof}

\section{The proofs of Theorems \ref{theorem1.1} and \ref{theorem1.2}} The proof of our 
main theorems using Theorem \ref{theorem6.2} is relatively routine, involving an 
optimisation of parameters. We first compute the Weyl-type exponent $\tau(k)$ defined in 
\eqref{5.2}. This is essentially the optimisation performed in the proofs of Corollaries 1 and 
2 to \cite[Theorem 1.1]{Woo1993}.\par

We begin by observing that whenever $v$ is even, then the corollary to 
\cite[Theorem 2.1]{Woo1993} shows that the exponent $\Del_v$ is admissible for $k\ge 4$, 
where $\Del_v$ is the unique positive solution of the equation
\begin{equation}\label{7.1}
\Del_ve^{\Del_v/k}=ke^{1-v/k}.
\end{equation}
Notice here that the exponent $\Del_s$ in the statement of this earlier result corresponds 
to our $\Del_v$ with $v=2s$, owing to the slightly different definitions employed between 
\cite{Woo1993} and the present memoir. Equipped with these exponents, we now seek to 
obtain a good approximation to
\begin{equation}\label{7.2}
k\tau(k)=\max_{w\in \dbN}\frac{1-2\Del_{2w}/k}{4w^2/k^2}.
\end{equation}
We explore this quantity by putting $w=\lceil \gam k\rceil$, where $\gam>0$ is a real 
parameter at our disposal. With the relation \eqref{7.1} in mind, we take $\del=\del(\gam)$ 
to be the positive solution of the equation
\begin{equation}\label{7.3}
\del+\log \del =1-2\gam .
\end{equation}
We note that the function $t+\log t$ is increasing for $t>0$. Then, since the relation 
\eqref{7.1} shows that the exponent $\Del_{2w}$ is admissible, where $\Del_{2w}$ is the 
unique positive solution of the equation
\[
\frac{\Del_{2w}}{k}+\log \frac{\Del_{2w}}{k}=1-\frac{2w}{k},
\]
and $1-2w/k\le 1-2\gam$, we infer that $\Del_{2w}\le k\del(\gam)$. We now define 
$\tet=\tet(\gam,w)$ by setting $\tet=w-\gam k$. Thus $0\le \tet<1$, and we see that the 
formula \eqref{7.2} delivers the lower bound
\begin{equation}\label{7.4}
k\tau(k)\ge \max_{\gam>0}\frac{1-2\del (\gam)}{4(\gam+\tet/k)^2}.
\end{equation}

\par One may now attempt to optimise the choice of $\gam$ on the right hand side of 
\eqref{7.4} so as to maximise our lower bound for $\tau(k)$. It transpires that the optimal 
choice for $\gam$ is very close to $1$, and so a good approximation to the maximum is 
found by taking $\gam=1$ and hence $\tet=0$. Solving \eqref{7.3} with $\gam=1$, it is 
apparent that $\del$ is constrained to satisfy the equation
\[
\del+\log \del +1=0.
\]
It is not difficult via a Newton iteration to verify that $\del=0.2784645\ldots $. With this 
value of $\del$, one has
\begin{equation}\label{7.5}
k\tau(k)\ge \frac{1-2\del}{4}=\frac{1}{9.027900\ldots }.
\end{equation}

\par Asymptotic information very slightly superior to the lower bound \eqref{7.5} is obtained 
by observing that since $0\le \tet<1$, the relation \eqref{7.4} yields
\[
k\tau(k)\ge \max_{\gam>0}\frac{1-2\del(\gam)}{4(\gam+1/k)^2}.
\]
The maximum here corresponds to a value of $\gam$ for which
\[
\frac{4(\gam+1/k)^2}{1-2\del(\gam)}
\]
achieves its minimum. On making use of \eqref{7.3} to eliminate $\gam$ and substituting 
$\xi$ for $\del(\gam)$, we find that this minimum value is equal to the minimum of the 
function
\[
\kap (\xi)=\frac{(1-\xi-\log \xi +2/k)^2}{1-2\xi},
\]
as $\xi$ varies over the interval $(0,1)$, and that the minimising value of $\del(\gam)$ is 
then equal to the value of $\xi$ corresponding to this minimum. Identifying the value of 
$\xi$ where $\kap'(\xi)=0$, we see that $\xi$ satisfies the equation
\[
\xi-\frac{1}{\xi}+2+\frac{2}{k}=\log \xi .
\]
Thus, if $\ome=3.548292\ldots $ is the positive real number with $\ome\ge 1$ satisfying 
the equation \eqref{1.1}, namely $\ome -2-1/\ome=\log \ome$, then we find that 
$\xi=1/\ome +O(1/k)$. We should therefore take $\del$ asymptotically close to $1/\ome$ 
for large $k$.\par

Motivated by this discussion, we put
\[
\gam=\tfrac{1}{2}(1-1/\ome+\log \ome)=0.992320\ldots ,
\]
and we avoid adjusting this value by the term of size $O(1/k)$ corresponding to the 
optimal choice. With this very slightly non-optimal choice of $\gam$, we find that
\[
k\tau(k)\ge \frac{1-2\del(\gam)}{4(\gam+1/k)^2}.
\]
Here, in view of \eqref{7.3}, one has
\[
\del+\log \del=1-2\gam=\frac{1}{\ome}+\log \frac{1}{\ome},
\]
whence $\del=1/\ome$. Thus
\[
k\tau(k)\ge \frac{1-2/\ome}{(1-1/\ome +\log \ome +2/k)^2}
=\frac{1}{9.026725\ldots }+O\Bigl(\frac{1}{k}\Bigr).
\]
We summarise these deliberations in the form of a lemma.

\begin{lemma}\label{lemma7.1}
When $k\ge 4$, one has
\[
\tau(k)\ge \frac{1}{9.027901 k},
\]
and also
\[
\tau(k)\ge \frac{1-2/\ome}{(1-1/\ome+\log \ome+2/k)^2k},
\]
where $\ome$ is the unique real solution with $\ome\ge 1$ of the equation
\[
\ome-2-1/\ome=\log \ome .
\]
\end{lemma}

We remark that, following a modest computation, one can confirm that the second lower 
bound for $\tau(k)$ delivered by this lemma takes the asymptotic form
\[
\tau(k)\ge \frac{1}{(\ome^2-3-2/\ome)k}+O\Bigl( \frac{1}{k^2}\Bigr) .
\] 

We may now make use of Theorem \ref{theorem6.2}, where we must consider the 
quantity
\[
G_0(k)=\min_{v\ge 2}\biggl( v+\frac{\Del_v}{\tau(k)}\biggr) .
\] 
Write $\tau(k)=(Dk)^{-1}$, where $D$ may depend on $k$, but is asymptotic to a 
constant determined via the conclusion of Lemma \ref{lemma7.1}. Then, on applying the 
formula \eqref{7.1} along the lines delivering \eqref{7.2}, we see that when $v$ is even 
one has
\begin{equation}\label{7.6}
v+\frac{\Del_v}{\tau(k)}\le v+Dk^2e^{1-\del -v/k},
\end{equation}
where $\del+\log \del=1-v/k$. As a corresponding inequality in a real variable $v$, the right 
hand side is approximately minimised by taking $v=k(1+\log (Dk))$. Instead, with 
$v$ constrained to be an even integer, we take
\[
v=2\left\lfloor \frac{1}{2}k(1+\log (Dk))-\frac{1}{2D}\right\rfloor .
\]
In this way, one finds that
\[
\del +\log \del \ge 1-(1+\log (Dk))+\frac{1}{Dk}=\frac{1}{Dk}+
\log \Bigl(\frac{1}{Dk}\Bigr),
\]
whence $\del\ge 1/(Dk)$.\par

Define the real number $\tet$ via the relation
\[
v=k(1+\log (Dk))-\frac{1}{D}-\tet ,
\]
and note that one then has $0\le \tet<2$. In this way, we discern that
\begin{align*}
1-\del-\frac{v}{k}&\le 1-\frac{1}{Dk}-(1+\log (Dk))+\frac{1}{Dk}+\frac{\tet}{k}\\
&=-\log (Dk)+\frac{\tet}{k}.
\end{align*}
Then we deduce from \eqref{7.6} that
\begin{equation}\label{7.7}
v+\frac{\Del_v}{\tau(k)}\le k(1+\log (Dk))-\frac{1}{D}-\tet+ke^{\tet/k}.
\end{equation}
The function $-\tet+ke^{\tet/k}$ is increasing with $\tet$ for $\tet\in [0,2)$, so is bounded 
above in this interval by $-2+ke^{2/k}$. Moreover, the function $k(e^{2/k}-1)$ is 
decreasing as a function of $k$ for $k\ge 2$. One may check that when $k\ge 20$, one has
\[
-2+ke^{2/k}\le k+\frac{1}{9.6694}<k+\frac{1}{D}.
\]
In such circumstances, we deduce that
\[
-\frac{1}{D}-\tet+ke^{\tet/k}\le -\frac{1}{D}-2+ke^{2/k}<k,
\]
whence, as a consequence of \eqref{7.7}, we obtain the bound
\[
v+\frac{\Del_v}{\tau(k)}<k(2+\log (Dk)).
\]
In this way, we deduce that for $k\ge 20$, one has
\begin{equation}\label{7.8}
G_0(k)\le k(\log k+2+\log D).
\end{equation}

\begin{proof}[The proof of Theorem \ref{theorem1.1}] By reference to the first bound 
supplied by Lemma \ref{lemma7.1}, one finds that the argument just described may be 
applied with $D=9.027901$ whenever $k\ge 20$. In such circumstances, one has 
$2+\log D\le 4.2003199$, and hence it follows from \eqref{7.8} that 
\[
\lfloor G_0(k)\rfloor \le \lfloor k(\log k+4.2003199)\rfloor \le \lceil k(\log k+4.20032)\rceil -1.
\]
The proof of Theorem \ref{theorem1.1} when 
$k\ge 20$ is therefore made complete by reference to Theorem \ref{theorem6.2}. For 
small values of $k$, one finds that the bounds for $G(k)$ already available in the literature 
are smaller than $\lceil k(\log k+4.20032)\rceil$ for $k\le 19$. Indeed, the bound 
$G(k)\le 2^k+1$ available via Hua's work (see the corollary to 
\cite[Theorem 2.6]{Vau1997}, for example) already suffices for $k\le 4$, while for 
$k\ge 14$ one has the bounds already reported in the introduction following the 
announcement of Theorem \ref{theorem1.3}. We can complete this list with the addition of 
the bounds $G(7)\le 31$, $G(8)\le 39$, $G(9)\le 47$, $G(10)\le 55$, $G(11)\le 63$, 
$G(12)\le 72$, $G(13)\le 81$, available from \cite{Woo2016}, together with the bounds 
$G(5)\le 17$ and $G(6)\le 24$ obtained, respectively, in \cite{VW1995} and \cite{VW1994}. 
Following this small list of checks, the proof of Theorem \ref{theorem1.1} is complete. 
\end{proof}

We note that the bound supplied by Theorem \ref{theorem1.1} is surprisingly competitive 
even for small values of $k$. Thus, for example, the bound $G(20)\le 144$ of Theorem 
\ref{theorem1.1} may be compared with the corresponding bound $G(20)\le 142$ of 
\cite{VW2000}. Of course, in Theorem \ref{theorem1.3} of the present memoir, we obtain 
$G(20)\le 137$.\par

\begin{proof}[The proof of Theorem \ref{theorem1.2}] We now apply the second bound 
supplied by Lemma \ref{lemma7.1}. With this bound in hand, the argument leading to 
\eqref{7.8} may be applied with
\[
D=\frac{(\ome-1-2/\ome+2/k)^2}{1-2/\ome},
\]
again, whenever $k\ge 20$. On recalling the definition \eqref{1.2} of $C_1$ and $C_2$, we 
now have
\begin{align*}
2+\log D&=2+\log \Bigl( \ome^2-3-\frac{2}{\ome}\Bigr) +
2\log \Bigl( 1+\frac{2}{k(\ome -1-2/\ome)}\Bigr) \\
&< C_1+\frac{4\ome}{k(\ome^2-\ome -2)}=C_1+\frac{C_2-1}{k}.
\end{align*}
We therefore deduce from \eqref{7.8} that
\[
G_0(k)+1<k(\log k+C_1+(C_2-1)/k)+1=k(\log k+C_1)+C_2.
\]
The proof of Theorem \ref{theorem1.2} is completed by reference to Theorem 
\ref{theorem6.2} when $k\ge 20$. For the small values of $k$ with $k\le 19$, the bound 
claimed in the statement of Theorem \ref{theorem1.2} is again confirmed by reference to 
the previously known upper bounds for $G(k)$ already cited in the proof of Theorem 
\ref{theorem1.1}.
\end{proof}

\section{Bounding $G(k)$ for intermediate values of $k$} Our proof of Theorem 
\ref{theorem1.3} follows the argument used to establish Theorems \ref{theorem1.1} and 
\ref{theorem1.2}, save that we now make use of the numerical tables of exponents 
available from \cite{VW2000}. We begin by numerically computing the exponent $\tau(k)$.

\begin{theorem}\label{theorem8.1} When $14\le k\le 20$, one has $\tau(k)\le T(k)^{-1}$, 
where the exponents $T(k)$ are presented in Table \ref{tab2}.
\end{theorem}

\begin{proof} We apply the formula
\[
T(k)=\Bigl( \frac{k-2\Del_{2w}}{4w^2}\Bigr)^{-1},
\]
available from \eqref{5.2}, using the values of $w$ and corresponding admissible exponents 
$\Del_{2w}$ to be found in the tables of \cite{VW2000}. Here, the exponents $\lam_w$ of 
\cite{VW2000} are related to $\Del_{2w}$ via the formula $\Del_{2w}=\lam_w-2w+k$. We 
record the necessary choice of parameter $w$, together with the associated admissible 
exponent $\Del_{2w}$, rounded up in the final decimal place, in Table \ref{tab2} below.
\end{proof}

\begin{table}[h]
\begin{center}
\begin{tabular}{ccccccc}
\toprule
$k$ & $w$  & $\Delta_{2w}$ & $T(k)$ & $v$ & $\Del_v$ & $G_0(k)$  \\
\toprule
$14$ & $26$ & $4.039939$ & $114.1869$ & $76$ & $0.109356$ & $88.4871$\\
$15$ & $28$ & $4.323087$ & $123.3903$ & $82$ & $0.117123$ & $96.4519$\\
$16$ & $30$ & $4.606286$ & $132.5981$ & $90$ & $0.108806$ & $104.4275$\\
$17$ & $32$ & $4.888677$ & $141.7763$ & $96$ & $0.116203$ & $112.4749$\\
$18$ & $34$ & $5.170691$ & $150.9411$ & $104$ & $0.109619$ & $120.5461$\\
$19$ & $36$ & $5.451758$ & $160.0695$ & $110$ & $0.116770$ & $128.6914$\\
$20$ & $38$ & $5.732224$ & $169.1748$ & $118$ & $0.111388$ & $136.8441$\\
\bottomrule
\end{tabular}\\[6pt]
\end{center}
\caption{Choice of exponents for $13\le k\le 20$.}\label{tab2}
\end{table}

We next confirm Theorem \ref{theorem1.3} by utilising the formula 
$G(k)\le \lfloor G_0(k)\rfloor +1$ available via Theorem \ref{theorem6.2}. Here, we have
\[
G_0(k)=v+\frac{\Del_v}{\tau(k)}=v+T(k)\Del_v,
\]
for a suitably chosen value of $v$. We present values of $v$, $\Del_v$ and $G_0(k)$ in 
Table \ref{tab2}, with the values $\Del_v$ extracted from \cite{VW2000}, again all rounded 
up in the final decimal place presented. The conclusion of Theorem \ref{theorem1.3} 
follows on noting that $G(k)\le \lfloor G_0(k)\rfloor +1$ for each value of $k$ in the table. 
This completes the proof of Theorem \ref{theorem1.3}.

\section{Remarks on upper bounds for $G^+(k)$} Scholars of the circle method as it 
applies to Waring's problem will appreciate instantly that the methods of this paper deliver 
bounds for the number $G^+(k)$, the smallest number $s$ having the property that almost 
all positive integers (in the sense of natural density) are the sum of at most $s$ positive 
integral $k$-th powers. Here, one makes a standard application of Bessel's inequality to 
estimate the minor arc contribution in mean square, the upshot being the familiar upper 
bound $G^+(k)\le \tfrac{1}{2}(H(k)+1)$, whenever $H(k)$ is an upper bound for $G(k)$ 
obtained by the methods of this paper. The methods here have nothing to contribute to the 
literature well-known to any worker in the area, so we may record without further delay the 
following conclusions.

\begin{theorem}\label{theorem9.1} Suppose that $k\in \dbN\setminus \{4,8,16,32\}$. Then
\[
G^+(k)\le \lceil \tfrac{1}{2}k(\log k+4.20032)\rceil
\]
and
\[
G^+(k)<\tfrac{1}{2}k(\log k+C_1)+\tfrac{1}{2}(C_2+1).
\]
In the exceptional cases $k=2^j$ with $j\in \{2,3,4,5\}$, one has $G^+(k)=4k$. Moreover, 
when $14\le k\le 20$ but $k\ne 16$, one has $G^+(k)\le H^+(k)$, where $H^+(k)$ is 
defined by means of Table \ref{tab3}.
\end{theorem}

\begin{table}[h]
\begin{center}
\begin{tabular}{ccccccccccccccccc}
\toprule
$k$ & $14$ & $15$ & $16$ & $17$ & $18$ & $19$ & $20$\\
$H^+(k)$ & $45$ & $49$ & $53$ & $57$ & $61$ & $65$ & $69$\\
\bottomrule
\end{tabular}\\[6pt]
\end{center}
\caption{Upper bounds for $G^+(k)$ when $14\le k\le 20$.}\label{tab3}
\end{table}

The assertion that $G^+(k)=4k$ when $k=2^j$ with $j\in \{2,3,4,5\}$ is not new. This was 
established by Hardy and Littlewood \cite{HL1925} when $k=4$, by Vaughan 
\cite{Vau1986} when $k=8$, and by the second author \cite{Woo1992} when $k=16$ and 
$k=32$. It is straightforward, however, to establish the following refinements that more 
fully reflect the entry $H^+(16)=53$ from Table \ref{tab3}, and the upper bound implicitly 
obtained for $k=32$ in Theorem \ref{theorem9.1}.

\begin{theorem}\label{theorem9.2} Let $k$ be either $16$ or $32$, and put 
$H^+(16)=53$ and $H^+(32)=123$. Suppose that $s\ge H^+(k)$ and that $r$ is an 
integer with $1\le r\le s$. Then almost all positive integers $n$ with $n\equiv r\mmod{4k}$ 
are the sum of $s$ positive integral $k$-th powers.
\end{theorem}

The proof of this conclusion is once again routine for scholars of the circle method, and we 
refer the reader to earlier literature such as \cite{Vau1986} or \cite{Woo1992} for the 
ideas necessary to complete this exercise.

\bibliographystyle{amsbracket}
\providecommand{\bysame}{\leavevmode\hbox to3em{\hrulefill}\thinspace}

\end{document}